\documentclass[12pt]{amsart}
\usepackage{amsmath,amssymb,amsbsy,amsfonts,latexsym,amsopn,amstext,
                                               amsxtra,euscript,amscd}

\usepackage{url}
\usepackage[colorlinks,linkcolor=blue,anchorcolor=blue,citecolor=blue]{hyperref}
\usepackage{color}

\newfont{\teneufm}{eufm10}
\newfont{\seveneufm}{eufm7}
\newfont{\fiveeufm}{eufm5}
\newfam\eufmfam
                \textfont\eufmfam=\teneufm \scriptfont\eufmfam=\seveneufm
                \scriptscriptfont\eufmfam=\fiveeufm
\def\bbbc{{\mathchoice {\setbox0=\hbox{$\displaystyle\rm C$}\hbox{\hbox
to0pt{\kern0.4\wd0\vrule height0.9\ht0\hss}\box0}}
{\setbox0=\hbox{$\textstyle\rm C$}\hbox{\hbox
to0pt{\kern0.4\wd0\vrule height0.9\ht0\hss}\box0}}
{\setbox0=\hbox{$\scriptstyle\rm C$}\hbox{\hbox
to0pt{\kern0.4\wd0\vrule height0.9\ht0\hss}\box0}}
{\setbox0=\hbox{$\scriptscriptstyle\rm C$}\hbox{\hbox
to0pt{\kern0.4\wd0\vrule height0.9\ht0\hss}\box0}}}}
\def\bbbq{{\mathchoice {\setbox0=\hbox{$\displaystyle\rm
Q$}\hbox{\raise 0.15\ht0\hbox to0pt{\kern0.4\wd0\vrule
height0.8\ht0\hss}\box0}} {\setbox0=\hbox{$\textstyle\rm
Q$}\hbox{\raise 0.15\ht0\hbox to0pt{\kern0.4\wd0\vrule
height0.8\ht0\hss}\box0}} {\setbox0=\hbox{$\scriptstyle\rm
Q$}\hbox{\raise 0.15\ht0\hbox to0pt{\kern0.4\wd0\vrule
height0.7\ht0\hss}\box0}} {\setbox0=\hbox{$\scriptscriptstyle\rm
Q$}\hbox{\raise 0.15\ht0\hbox to0pt{\kern0.4\wd0\vrule
height0.7\ht0\hss}\box0}}}}
\def\bbbt{{\mathchoice {\setbox0=\hbox{$\displaystyle\rm
T$}\hbox{\hbox to0pt{\kern0.3\wd0\vrule height0.9\ht0\hss}\box0}}
{\setbox0=\hbox{$\textstyle\rm T$}\hbox{\hbox
to0pt{\kern0.3\wd0\vrule height0.9\ht0\hss}\box0}}
{\setbox0=\hbox{$\scriptstyle\rm T$}\hbox{\hbox
to0pt{\kern0.3\wd0\vrule height0.9\ht0\hss}\box0}}
{\setbox0=\hbox{$\scriptscriptstyle\rm T$}\hbox{\hbox
to0pt{\kern0.3\wd0\vrule height0.9\ht0\hss}\box0}}}}
\def\bbbs{{\mathchoice
{\setbox0=\hbox{$\displaystyle     \rm S$}\hbox{\raise0.5\ht0\hbox
to0pt{\kern0.35\wd0\vrule height0.45\ht0\hss}\hbox
to0pt{\kern0.55\wd0\vrule height0.5\ht0\hss}\box0}}
{\setbox0=\hbox{$\textstyle        \rm S$}\hbox{\raise0.5\ht0\hbox
to0pt{\kern0.35\wd0\vrule height0.45\ht0\hss}\hbox
to0pt{\kern0.55\wd0\vrule height0.5\ht0\hss}\box0}}
{\setbox0=\hbox{$\scriptstyle      \rm S$}\hbox{\raise0.5\ht0\hbox
to0pt{\kern0.35\wd0\vrule height0.45\ht0\hss}\raise0.05\ht0\hbox
to0pt{\kern0.5\wd0\vrule height0.45\ht0\hss}\box0}}
{\setbox0=\hbox{$\scriptscriptstyle\rm S$}\hbox{\raise0.5\ht0\hbox
to0pt{\kern0.4\wd0\vrule height0.45\ht0\hss}\raise0.05\ht0\hbox
to0pt{\kern0.55\wd0\vrule height0.45\ht0\hss}\box0}}}}
\def\bbbz{{\mathchoice {\hbox{$\sf\textstyle Z\kern-0.4em Z$}}
{\hbox{$\sf\textstyle Z\kern-0.4em Z$}} {\hbox{$\sf\scriptstyle
Z\kern-0.3em Z$}} {\hbox{$\sf\scriptscriptstyle Z\kern-0.2em
Z$}}}}

\newtheorem{theorem}{Theorem}

\newtheorem{prob}[theorem]{Problem}
\newtheorem{fact}[theorem]{Fact}

 \numberwithin{equation}{section}
  \numberwithin{theorem}{section}

\def\squareforqed{\hbox{\rlap{$\sqcap$}$\sqcup$}}
\def\qed{\ifmmode\squareforqed\else{\unskip\nobreak\hfil
\penalty50\hskip1em\null\nobreak\hfil\squareforqed
\parfillskip=0pt\finalhyphendemerits=0\endgraf}\fi}

\def\cI{{\mathcal I}}

\def\cL{{\mathcal L}}

\def\cP{{\mathcal P}}

\def \sf {\mathfrak s}

\def\fG {\mathfrak G}

\def\fP {\mathfrak P}

\newcommand{\ignore}[1]{}

\def\vec#1{\mathbf{#1}}

\def \Z{\mathbb{Z}}

\def \Z{\mathbb{Z}}

\def\mand{\qquad\mbox{and}\qquad}

\def\\{\cr}
\def\({\left(}
\def\){\right)}

\global\long\def\mF{\mathbb{F}}

\global\long\def\mZ{\mathbb{Z}}

\global\long\def\cI{\mathcal{I}}

\global\long\def\cL{\mathcal{L}}

\global\long\def\cP{\mathcal{P}}

\renewcommand{\vec}[1]{\mathbf{#1}}

\begin{document}

\title{Counting Co-Cyclic Lattices}

\author[P. Q. Nguyen] {Phong Q. Nguyen}
\address{
INRIA, Domaine de Voluceau, 78153 Rocquencourt, France and Tsinghua University, Institute for Advanced Study, Beijing 100084, China}
\email{phong.nguyen@inria.fr}
 
  \author[I. E. Shparlinski] {Igor E. Shparlinski}
\address{Department of Pure Mathematics, University of New South Wales,
Sydney, NSW 2052, Australia}
\email{igor.shparlinski@unsw.edu.au}

\begin{abstract} 
There is a well-known asymptotic formula, due to W.~M.~Schmidt (1968) for the number of full-rank integer 
lattices of index at most $V$ in $\mZ^n$.
This set of lattices $L$ can naturally be partitioned with respect to the factor  group $\mZ^n/L$.
Accordingly, we count the number of full-rank integer lattices $L \subseteq \mZ^n$ such that $\mZ^n/L$ is cyclic and of order at most $V$,
and deduce that these co-cyclic lattices are dominant among all integer lattices:
their natural density is  $\(\zeta(6) \prod_{k=4}^n \zeta(k)\)^{-1} \approx 85\%$.
The problem is motivated by complexity theory, namely worst-case to average-case reductions for lattice problems.
\end{abstract}


\maketitle

\section{Introduction}

Let $\cI_{n,V}$ be the set of  subgroups $L$ of $\mZ^n$  such that $[\mZ^n:L] = V$.
In other words, $\cI_{n,V}$ is the (finite) set of full-rank integer lattices $\subseteq \mZ^n$
of determinant (or co-volume) equal to $V$. Let 
$$
\cI_{n,\le V}= \cup_{1 \le v \le V} \cI_{n,v}.
$$
A classical result (see~\cite{GiGr06,Sc68}) states that when $n$ is fixed and $V$ grows to $\infty$,
\begin{align} \label{eq:nblattices}
  \#\cI_{n,\le V}  \sim \Xi_{2,n} \frac{V^n}{n}  
\end{align}
where we define
$$
\Xi_{m,n} = \prod_{k=m}^n \zeta(k) \mand
\Xi_{m} =  \prod_{k=m}^\infty \zeta(k)
$$
and $\zeta$ is the  {\it Riemann zeta-function\/}.
This is a special case of subgroup growth for the group $\mZ^n$.

We are interested in counting special subsets of $\cI_{n,\le V}$.
More precisely, $\cI_{n,\le V}$ can naturally be partitioned with respect to the (finite Abelian) factor group $\mZ^n/L$ for $L \in \cI_{n,\le V}$.
For any finite Abelian group $G$, we 
denote by $\cL_{n,G}$ the finite set of full-rank integer lattices $L \subseteq \mZ^n$
such that $\mZ^n/L \simeq G$. Then 
$$
\cI_{n,V} = \cup_{\#G = V} \cL_{n,G}\mand \cI_{n,\le V} = \cup_{\#G \le V} \cL_{n,G}.
$$
The sets $\cL_{n,G}$  have attracted significant interest in complexity theory.
In a seminal work~\cite{Aj96}, Ajtai has discovered the first worst-case to average-case reduction for lattice problems:
 he proves that for the special case $G = (\mZ/q\mZ)^m$, when $\#G$ is sufficiently large,
finding very short non-zero vectors (with non-negligible probability) in a random lattice $L \in \cL_{n,G}$ chosen with uniform distribution
is as hard as finding short non-zero vectors in any lattice of dimension $m$.
Ajtai's reduction has recently been generalised~\cite{GIN} to any finite Abelian group $G$ of sufficiently large order,
which motivates to study the cardinals of $\cL_{n,G}$, depending on $G$.

We settle this question for cyclic groups $G$.
More precisely,  we give an asymptotic formula for the cardinality 
$N_n(V)$ of the the subset of $\cI_{n,\le V}$ formed by all {\em co-cyclic lattices}, 
that is, full-rank integer lattices $L$ such that $\mZ^n/L$ is cyclic:
$$N_n(V) = \# \cup_{\substack{G~\mathrm{cyclic}\\ \#G \le V}} \cL_{n,G},
$$
which is the subset of $\cI_{n,\le V}$ formed by all {\em co-cyclic lattices}, that is, full-rank integer lattices $L$ such that $\mZ^n/L$ is cyclic.
Such lattices have been previously studied from a complexity point 
of view in~\cite{PaSch,Trol}.

Throughout the paper we use the notion of {\em natural density\/}. 
We recall that the natural density of a property $\cP$ in a family of objects  ordered according to their ``size'' (such as lattices of  
determinant up to $V$, 
groups of order up to $V$ and so on) is defined as the limit as $V\to \infty$
of the ratio of the number of such objects of size at most $V$ satisfying the property $\cP$
 to the total number of such objects of size at most $V$.

For example, our results, coupled with~\eqref{eq:nblattices},
show that the natural density of co-cyclic lattices of  a fixed dimension $n$ 
 tends  asymptotically (as $n$  grows) to 
\begin{align}
\label{eq:density}
\frac{1}{\zeta(6)  \Xi_{4} } \approx 85\%. 
\end{align}
Hence, ``most'' integer lattices of sufficiently large dimension are co-cyclic.

We also determine the density of integer lattices with squarefree index in $\mZ^n$,
which are special cases of co-cyclic lattices, 
see also the discussion in Section~\ref{sec:gen}.
More precisely, we obtain an asymptotic formula on the
cardinality
$$
N^\sharp_n(V) = \# \cup_{\substack{\#G \le V\\\#G~\mathrm{squarefree}}} \cL_{n,G}.
$$
Coupled with~\eqref{eq:nblattices},
we obtain that the natural density of full-rank integer lattices  a fixed dimension $n$  
with squarefree determinant of  a fixed dimension $n$
 tends  asymptotically (as $n$  grows) to  
\begin{align}
\label{eq:density sf}
\frac{1}{\Xi_{3}} \approx 71.7\%.
\end{align}

\section{Main Results}

For $n \ge 2$ we define $\vartheta_n$ by the absolutely converging product
$$
\vartheta_n  
=  \prod_{p} \(1 + \frac{p^{n-1} -1}{p^{n+1} - p^{n}}\), 
$$
where hereafter $p$ always runs through prime numbers.

\begin{theorem}
\label{thm:NnV}
For any fixed $n \ge 2$, we have
$$
N_n(V) = \frac{\vartheta_n }{n} V^n   +
O\( V^{n-1+o(1)}\),
$$
as $V\to \infty$.
\end{theorem}

Similarly, we define
$$
\rho_n  
=  \frac{6}{\pi^2} \prod_{p} \(1 + \frac{p^{n-1} -1}{p^{n+1} - p^{n-1}}\) . 
$$

\begin{theorem}
\label{thm:NnV-SF}
For any fixed $n \ge 2$, we have
$$
N^\sharp_n(V) = \frac{\rho_n +o(1)}{n} V^n 
$$
as $V\to \infty$.
\end{theorem}

We note that it is easy to get an explicit bound on error term 
in the asymptotic formula of Theorem~\ref{thm:NnV-SF}.

It is also interesting to study the asymptotic behaviour of the constant $\vartheta_n$ as
$n \to \infty$. We define
$$ \vartheta = \prod_p \left( 1+ \frac{1}{p^2-p} \right) = \frac{\zeta(2)\zeta(3)}{\zeta(6)}= \frac{ 315 \zeta(3)}{2 \pi^4}  =1.94359\ldots$$
We note that $\vartheta$ appears in  the asymptotic formula of
Landau~\cite{Lan}:
\begin{equation}
\label{eq:Landau}
\sum_{d\le t} \frac{1}{\varphi(d)} =
\vartheta \left( \log t + \gamma -
\sum_{p}\frac{\log p}{p^2-p+1} \right ) + O\( \frac{\log t}{t} \),
\end{equation}
where $\varphi(d)$ is the Euler function and 
$\gamma$ is the {\it Euler-Mascheroni constant\/} (a more recent reference
is~\cite{Mont}).

\begin{theorem}
\label{thm:theta}
For any  $n \ge 2$, we have
$$
\vartheta\(1 -   \frac{7\cdot 2^{-n}}{6} +O(3^{-n})\) \ge \vartheta_n \ge \vartheta\(1 -  
 \frac{5\cdot 2^{-n}}{3}  +O(3^{-n})\)  . 
$$
\end{theorem}
By combining Theorems~\ref{thm:NnV} and~\ref{thm:theta} with~\eqref{eq:nblattices},
we obtain~\eqref{eq:density}.

We now define 
$$ \rho= \prod_p \left( 1+ \frac{1}{p^2-1} \right) = \zeta(2).
$$
Finally, we also have

\begin{theorem}
\label{thm:rho}
For any  $n \ge 2$, we have
$$
\rho_n = \rho\(1 -    2^{-n-1} +O(3^{-n})\)  . 
$$
\end{theorem}

By combining Theorems~\ref{thm:NnV} and~\ref{thm:rho} with~\eqref{eq:nblattices},
we obtain~\eqref{eq:density sf}.

One can easily get tighter bounds in Theorems~\ref{thm:theta} and~\ref{thm:rho}.

\section{Proofs of Main Results}

\subsection{Proof of Theorem~\ref{thm:NnV}}
\label{sec:proof1}

Given an integer $q$, we say that two $n$-dimensional 
vectors $\vec{a}, \vec{b} \in \Z^n$ are {\it equivalent\/} modulo $q$,
if for some integer  $\lambda$ with $\gcd(\lambda,q) =1$ we have
$$
\vec{a} \equiv \lambda \vec{b} \pmod q.
$$
We also say that a vector $\vec{a}=(a_1,\ldots, a_n) \in \Z^n$ 
is {\it primitive\/} modulo $q$,
if $\gcd(a_1,\ldots, a_n, q) = 1$. 

Let $A_n(q)$ be the number of distinct non-equivalent primitive modulo $q$
vectors $\vec{a}=(a_1,\ldots, a_n) \in \Z^n$. 

Let $L \in \cI_{n,q}$. Paz and Schnorr~\cite{PaSch} have proved that $L$ is co-cyclic if and only if there exist
a vector  $\vec{a}=(a_1,\ldots, a_n) \in \Z^n$   primitive modulo $q$,
such that $L = \cL_n(q,\vec{a})$ where
\begin{equation}
\begin{split}
\label{eq:Lnq}
 \cL_n(q, \vec{a})  = \{\vec{x} &=(x_1,\ldots, x_n) \in \Z^n~:\\
 &~a_1 x_1 + \ldots + a_nx_n  \equiv 0  \pmod q\}. 
\end{split}
\end{equation}
It follows that $N_n(V)$ satisfies:
$$
N_n(V) = \sum_{q \le V} A_n(q)
$$
Using the M{\"o}bius function $\mu(d)$, see~\cite[Section~16.3]{HW},
we write
\begin{equation}
\label{eq:An1}
A_n(q) =\frac{1}{\varphi(q)}  \sum_{d\mid q}\mu(d)
\sum_{\substack{a_1, \ldots, a_n =1 \\ 
 d \mid  \gcd(a_1,\ldots, a_n)}}^q 1
=\frac{q^n}{\varphi(q)}  \sum_{d\mid q}\frac{\mu(d)}{d^n}.
\end{equation}

Using the well-known identity (see~\cite[Theorem~62]{HW})
$$
\varphi(q) = q \prod_{p\mid q}\(1 - \frac{1}{p}\),
$$
and since
$$
\sum_{d\mid q}\frac{\mu(d)}{d^n} =  \prod_{p\mid q}\(1 - \frac{1}{p^n}\)
$$
we derive from~\eqref{eq:An1}
\begin{equation}
\label{eq:An2}
A_n(q)= q^{n-1}  \prod_{p\mid q}\(1 + \frac{p^{n-1} -1}{p^n - p^{n-1}}\).
\end{equation}

Let $f_n(d)$ be 
the 
multiplicative function  defined as 
$$
f_n(d) = \prod_{p \mid d} \frac{p^{n-1} -1}{p^n - p^{n-1}}
$$
if $d$ is squarefree and  $f(d) =  0$ otherwise. 
From~\eqref{eq:An2}, we see   that
$$
A_n(q)= q^{n-1} \sum_{d \mid q} f_n(d).
$$
Therefore, changing the order of summation and writing $q = kd$, we
derive
\begin{equation}
\label{eq:Nn}
N_n(V) = \sum_{d \le V} f_n(d) d^{n-1} \sum_{k \le V/d} k^{n-1}.
\end{equation}
We now observe that 
$$
\frac{p^{n-1} -1}{p^n - p^{n-1}} \le \frac{2}{p}.
$$
Hence, for any integer  $d$  we have 
$$
f(d) \le \frac{2^{\omega(d)}}{d}, 
$$
where $\omega(d)$ is the number of prime divisors of $d$.
Recalling the well-known bound
$$
\omega(d) = O\(\frac{\log d}{\log \log d}\)
$$
(which follows immediately from the trivial inequality $\omega(d)! \le d$ 
and Stirling's formula) we obtain
\begin{equation}
\label{eq:fd}
f(d) \le  d^{-1+o(1)}, \qquad \text{as}\ d\to \infty. 
\end{equation}

Since, for a fixed $n$, 
$$
\sum_{k \le V/d} k^{n-1} = \frac{1}{n} (V/d)^n + O\((V/d)^{n-1}\)
$$
we now derive from~\eqref{eq:Nn} and~\eqref{eq:fd} that
\begin{equation}
\label{eq:Nn prelim}
\begin{split}
N_n(V) &= \frac{1}{n} V^n \sum_{d \le V} \frac{f_n(d)}{d}  +
O\( V^{n-1} \sum_{d \le V} f_n(d)\)\\
& = \frac{1}{n} V^n \sum_{d \le V} \frac{f_n(d)}{d}  +
O\( V^{n-1+o(1)}\) ,
\end{split}
\end{equation}
as $V\to \infty$. Finally, using~\eqref{eq:fd} again, we obtain
$$
\sum_{d \le V} \frac{f_n(d)}{d} = \vartheta_n  + O\( V^{n-1+o(1)}\),
$$
where $\vartheta_n$ is given by the absolutely converging series
$$
\vartheta_n = \sum_{d =1}^\infty  \frac{f_n(d)}{d} = 
\prod_{p} \(1 +  \frac{f_n(p)}{p}\) 
=  \prod_{p} \(1 + \frac{p^{n-1} -1}{p^{n+1} - p^{n}}\).
$$
which concludes the proof. 
 
\subsection{Proof of Theorem~\ref{thm:NnV-SF}}

For two functions $F(t)$ and $G(t)$ depending on a parameter $t$ we write 
$F(t) \sim G(t)$ as an equivalent of 
$$
\lim_{t\to \infty} F(t)/G(t) = 1.
$$

We need the following classical result of Prachar~\cite{Prach} which asserts
that for any integers $d> a \ge 1$ with $\gcd(d,a)=1$ and arbitrary $\varepsilon >0$
we have
$$
\# \{k\le x~:~k\equiv a \pmod d, \ k~\text{squarefree}\}
 \sim \frac{6x}{\pi^2d}\prod_{p\mid d}\(1-\frac{1}{p^2}\)^{-1},
$$
as $x \to \infty$,
provided $d \le x^{2/3-\varepsilon}$, where the product is taken
over all prime divisors $p$ of $d$; see also~\cite[Theorem~3]{Hool}. 
In particular, under the same condition on $d$ and $x$, we have
\begin{equation}
\label{eq:SF coprime}
\begin{split}
\# \{k\le x~ :~&\gcd(k,d)=1,  \ k~\text{squarefree}\}\\
&\sim \frac{6\varphi(d)x}{\pi^2 d}\prod_{p\mid d}\(1-\frac{1}{p^2}\)^{-1}\sim \frac{6x}{\pi^2}\prod_{p\mid d}\(1+\frac{1}{p}\)^{-1},
\end{split}
\end{equation}
as $x \to \infty$ (one can certainly prove~\eqref{eq:SF coprime} directly 
as well).

As in the proof of  Theorem~\ref{thm:NnV} we write
$$ N^\sharp_n(V) = \sum_{\substack{q \le V\\ q~\text{squarefree}}} A_n(q) .
$$
Furthermore, instead of~\eqref{eq:Nn}, we derive
$$ N^\sharp_n(V) = \sum_{d \le V} f_n(d) d^{n-1}
\sum_{\substack{k \le V/d\\\gcd(k,d)=1\\ k~\text{squarefree}}}  k^{n-1}
$$
(recall that the function $f(d)$ is supported only on squarefree integers). 

For $d > V^{1/2}$ we estimate the sum over $k$ trivially as $O((V/d)^{n})$.
Thus, recalling~\eqref{eq:fd}, we see that the total contribution from 
such terms is  
\begin{equation}
\label{eq:large d}
\begin{split}
\sum_{V^{1/2} < d \le V} f_n(d) &d^{n-1} 
\sum_{\substack{k \le V/d\\\gcd(k,d)=1\\ k~\text{squarefree}}}  k^{n-1}\\
&= 
V^n \sum_{V^{1/2} < d \le V}  d^{-2+o(1)}
=  V^{n-1/2+o(1)}.
\end{split}
\end{equation}

For $d \le  V^{1/2}$, the asymptotic formula~\eqref{eq:SF coprime}
applies to the sums over $k$, so similarly to~\eqref{eq:Nn prelim},
via partial summation, 
we derive that the total contribution from 
such terms is  
\begin{equation}
\label{eq:small d}
\begin{split}
\sum_{d \le V^{1/2}} f_n(d) &d^{n-1} 
\sum_{\substack{k \le V/d\\\gcd(k,d)=1\\ k~\text{squarefree}}}  k^{n-1}\\
&\sim \frac{6  V^n }{\pi^2n }\sum_{d\le V^{1/2}} \frac{f_n(d)}{d}
\prod_{p\mid d}\(1+\frac{1}{p}\)^{-1},
\end{split}
\end{equation}
as $V\to \infty$. 

Since
\begin{equation*}
\begin{split}
\sum_{d\le V^{1/2}} \frac{f_n(d)}{d}
&\prod_{p\mid d}\(1+\frac{1}{p}\)^{-1}
  \sim \sum_{d=1}^\infty  \frac{f_n(d)}{d}
\prod_{p\mid d}\(1+\frac{1}{p}\)^{-1}  \\
& = \prod_{p} \(1 +  \frac{f_n(p)}{p} \(1+\frac{1}{p}\)^{-1} \)  \\
& = \prod_{p} \(1 +  \frac{f_n(p)}{p+1} \)  =  \prod_{p} \(1 + \frac{p^{n-1} -1}{p^{n+1} - p^{n-1}}\) , 
\end{split}
\end{equation*}
the result now follows from~\eqref{eq:large d} and~\eqref{eq:small d}. 

\subsection{Proof of Theorem~\ref{thm:theta}}

First we note that for any $k \ge 2$
\begin{equation}
\label{eq:zeta}
\zeta(k) = 1 +  \frac{1}{2^{k}} + \frac{1}{3^{k}} + O\(\int_{3}^\infty\frac{1}{t^{k}}dt\)
= 1 +  2^{-k} + O(3^{-k})  
\end{equation}
 We have,
$$ \frac{  \vartheta_n }{\vartheta} = \prod_p \frac{ p^2-p+1-1/p^{n-1}}{p^2-p+1} = \prod_p \( 1- \frac{1}{p^{n-1}(p^2-p+1)}\) .$$
Since  $p^2-p+1 \ge p$,  we have 
$$
1- \frac{1}{p^{n-1}(p^2-p+1)} \ge  1- \frac{1}{p^{n}} . 
$$
Hence, we now see that
\begin{equation*}
\begin{split}
 \frac{  \vartheta_n }{\vartheta} &\ge 
 \( 1- \frac{1}{3\cdot 2^{n-1}}\) \prod_{p\ge 3} \( 1- \frac{1}{p^{n}} \) \\
& = \( 1- \frac{1}{3\cdot 2^{n-1}}\) \( 1- \frac{1}{3^{n}}\)^{-1}  \zeta(n)^{-1}. 
\end{split}
\end{equation*}
Thus, using~\eqref{eq:zeta} we obtain 
$$
 \frac{  \vartheta_n }{\vartheta} \ge 
1- \frac{1}{3\cdot 2^{n-1}} - \frac{1}{2^n}   + O(3^{-n}) 
$$
On the other hand, since
$p^2-p+1 \le p^2$, we also have
\begin{equation*}
\begin{split}
 \frac{  \vartheta_n }{\vartheta} & \le   \( 1- \frac{1}{3\cdot 2^{n-1}}\) \prod_{p\ge 3}
 \( 1- \frac{1}{p^{n+1}}\)  \\
 & =  \( 1- \frac{1}{3\cdot 2^{n-1}}\) \( 1- \frac{1}{3^{n+1}}\)^{-1}\zeta(n+1)^{-1}.
\end{split}
\end{equation*}

\subsection{Proof of Theorem~\ref{thm:rho}}

 We have,
$$ \frac{\rho_n }{\rho} = \prod_p \( 1- \frac{1}{p^{n+1}}\) = \frac{1}{\zeta(n+1)}.
$$
Using~\eqref{eq:zeta}, we conclude the proof.

\section{Comparison with random finite Abelian groups}

\subsection{Motivation}
 In any context involving finite Abelian groups, it is interesting
to study if these finite Abelian groups behave like random finite Abelian groups
in various families in terms of natural density.

There are at least two distributions worth considering, which we discuss here.

 \subsection{The uniform distribution} 
Let $a(n)$ denote the number of non-isomorphic Abelian groups of order $n$.
It is well-known that if $n=p_1^{k_1} \ldots p_s^{k_s}$ is a prime number factorisation 
of $n$ then
$$
a(n) = \prod_{i=1}^s \cP(k_i)
$$
where $\cP(k)$ is the number integer partitions of $k$ (as it is obvious that $a(n)$ is a 
multiplicative function and also $a(p^k) = \cP(k)$ for any integer power of a prime $p$), and also 
$$  \sum_{\#G \le x} 1= \sum_{n \le x} a(n) = A_1 x + A_2 x^{1/2} + A_3 x^{1/3} + R(x),$$
where 
$$
A_j = \prod_{k=1, k \ne j}^{\infty} \zeta(k/j), \quad j=1,2,3,
$$
and the best result for the error term is $R(x) \ll x^{1/4+o(1)}$ 
(see~\cite{RoSa06}).

Since clearly there is a unique  (up to isomorphism) cyclic  subgroup of order $n$
and furthermore, all subgroups of square-free order are cyclic, we see that 
the natural density of  cyclic groups is 
$$
\frac{1}{A_1} =  \frac{1}{\Xi_2}  \approx 44\% .
$$
Furthermore, recalling that there are $\(\zeta(2)^{-1}+o(1)\)V$ 
squarefree number up to $V$ (see~\cite[Section~18.6]{HW}), we conclude 
that  the natural density of  groups of square-free order is 
$$
\frac{1}{\zeta(2) A_1} =\frac{1}{\zeta(2) \Xi_2}  \approx 26\% .
$$

\subsection{The Cohen-Lenstra distribution}  
According to the Cohen-Lenstra heuristics~\cite{CoLe84},
a given finite Abelian group $G$ occurs
  with mass inversely proportional
to the order  $\#\mathrm{Aut}(G)$ of its automorphism group $\mathrm{Aut}(G)$, similarly to  many other 
mass formulas.  
Let 
$$
T(V) = \sum_{\#G \le V} 1/\#\mathrm{Aut}(G),
$$
with 
$G$ running over all finite Abelian groups of order at most  $V$, up to isomorphism; 
Furthermore, let $T(V; \fP)$ is the same sum as $T(V)$ restricted to the groups $G$  satisfying $\fP$. 
Then the natural density of a property $\fP$ is 
 defined here as the limit  
\begin{equation}
\label{eq:Nat Dens}
\Delta(\fP) = \lim_{V\to \infty} T(V;  \fP)/T(V),
\end{equation}
provided it exists.

It is shown in~\cite{CoLe84} that the denominator is asymptotically equivalent to:
 $$
 T(V) \sim   \Xi_2 \log V.
 $$
If $\fP$ is the property $\fP_{cycl}=\text{``{\em $G$ is cyclic}''}$, then the numerator is:
\begin{align*}
T(V;  \fP_{cycl}) & = \sum_{n \le V} \frac{1}{\varphi(n)}\\
& =
\vartheta \left( \log V + \gamma -
\sum_{p}\frac{\log p}{p^2-p+1} \right ) + O\( \frac{\log x}{x} \),
\end{align*}
see~\eqref{eq:Landau}. 
Hence, with respect to the Cohen-Lenstra distribution,
the natural density of cyclic groups is 
\begin{equation}
\label{eq:cycl dens}
\Delta(\fP_{cycl}) = \vartheta  \Xi_2^{-1} = \zeta(6)^{-1} \Xi_4^{-1} \approx 85\%.
\end{equation}

Hence, the cyclicity of the factor group $\mZ^n/L$ (when $L$ is full-rank integer lattices $L$ in $\mZ^n$)
behaves as predicted by the Cohen-Lenstra heuristics.

If $\fP_{sf}$ is the property  $\fP_{sf}=\text{``{\em $\#G$ is square-free}''}$,
then:
$$
T(V;  \fP_{sf})  = \sum_{n \le V} \frac{\mu^2(n)}{\varphi(n)}.
$$
Ward~\cite[Equation~(2.2)]{Wa27} has shown that, as $V$ grows to $\infty$:
$$ \sum_{n \le V} \frac{\mu(n)^2}{\varphi(n)} =  \log V + c + o(V^{-1/2})
$$
for some absolute constant $c$. 
Hence, with respect to the Cohen-Lenstra distribution,
the natural density of groups of square-free order is 
$$ \Delta(\fP_{sf}) =  \Xi_2^{-1} \approx 44\%,$$
which differs from the natural density of co-cyclic lattices.

More generally, it is also possible to obtain the natural density when $\fP$ is the property
$\fP_{\le r}=\text{``{\em $G$ has rank at most $r$}''}$.
To do so, we rely on results from~\cite{CoLe84} for $p$-groups, for which the Cohen-Lenstra distribution is a probability distribution:
the probability that a random $p$-group 
has rank $r$ is exactly:
$$ P(p,r) = p^{-r^2} \frac{\prod_{i=1}^{\infty} (1-p^{-i})}{\prod_{i=1}^{r} (1-p^{-i})^2}.$$
It follows that the density of $\fP_{\le r}$  is:
\begin{align*}
 \Delta(\fP_{\le r}) & =  \prod_{p}   \sum_{k=0}^r P(p,k)\\ & =  \prod_{p}   
 \(  \sum_{k=0}^r p^{-k^2} \frac{(1-p^{-1})}
{\prod_{i=1}^{k} (1-p^{-i})^2}  \prod_{i=2}^{\infty} (1-p^{-i})  \)  \\
& =  \Xi_2^{-1} \prod_{p} \(    \sum_{k=0}^r p^{-k^2} \frac{(1-p^{-1})}
{\prod_{i=1}^{k} (1-p^{-i})^2} \).
\end{align*}
This proves again  that 
$\Delta(\fP_{\le 1}) =   \vartheta  \Xi_2^{-1}$,  
which is the previous natural density $\Delta(\fP_{cycl})$ of cyclic groups~\eqref{eq:cycl dens}.
We also have:
\begin{align*}
\Delta(\fP_{\le 2})& = \Xi_2^{-1} \prod_{p} \( 1+\frac{1}{p(p-1)} + \frac{p^{-4}}{(1-p^{-1})(1-p^{-2})^2}  \)  \\
 & =  \Xi_2^{-1}
  \prod_{p} \left( 1+ \frac{p^4-p^2+1}{p(p-1)^3(p+1)^2}   \right)  \approx 99.5\%.
 \end{align*}
 It is shown in~\cite{FuGo15} that for all $p \ge 2$:
 $$ 1-p^{-1}-p^{-2} \le \prod_{i\ge 1} (1-p^{-i}) \le 1.$$
It follows that for all $p \ge 2$:
$$ P(p,k) \le \frac{p^{-k^2}}{1- p^{-1}-p^{-2}} \le 4 p^{-k^2}.$$
Therefore, 
$$
\sum_{k\ge r} P(p,k) \le  \sum_{k \ge r} 4p^{-k^2} \le 8p^{-r^2}.
$$ 
(provided the above sums and products converge). Hence the natural density of the property $\fP_{\ge r}=\text{``{\em $G$ has rank at least $r$}''}$ satisfies:
\begin{align*}
 \Delta(\fP_{\ge r}) & = 1- \prod_{p}   \sum_{j=0}^{r-1} P(p,j)
 = 1- \prod_{p}  \(1 -  \sum_{j\ge r} P(p,j)\) \\
 & \le  1- \prod_{p}  \(1 -  8p^{-r^2}\) = 1- \exp\(\sum_{p} 8p^{-r^2}\)\\
& \le 1- \exp\(8 (\zeta(r^2)-1)\).
\end{align*}
Since $\zeta(r^2)-1) = 2^{-r^2} + O(3^{-r^2}))$, see~\eqref{eq:zeta}
we obtain 
$$ 
\Delta(\fP_{\ge r}) \le 8\cdot 2^{-r^2} + O(3^{-r^2}) .
$$

\subsection{Comparison with groups of points of elliptic curves over prime fields}

We note that Gekeler~\cite{Ge08} settled analogous questions for elliptic curves over prime finite fields:
if $q$ is a large random prime and $E$ runs over all elliptic curves over $\mF_q$,
then
\begin{itemize}
\item  The natural density of cyclic $E$ is  $$ \prod_{p} \left( 1 - \frac{1}{(p^2-1)p(p-1)} \right) \approx 0.81 .$$
\item The natural density of $E$ such that $\#E$ is squarefree is $$ \prod_{p} \left( 1 - \frac{p^3-p-1}{(p^2-1)p^2(p-1)} \right) \approx 0.44.$$
\end{itemize}

\section{Comments and Open Questions} 

\label{sec:gen}
We have studied the cardinality $N_n(V)$ of the 
 subset of $\cI_{n,\le V}$ formed by all lattices $L$ such that $\mZ^n/L$ is cyclic.
More generally, we define $\fG_m(V)$ as the set 
of groups $G$ with  $m$ invariant factors (that is, of rank $m$, where the rank is defined as the minimal size of a generating set) 
and of order $\# G \le V$.
It is natural to study the cardinality $N_{n,m}(V)$ of the 
subset $\cup_{G\in \fG_m(V)} \cL_{n,G}$ of $\cI_{n,V}$
 formed by all lattices $L$ such that $\mZ^n/L$ has exactly $m$ invariant factors.
 We are now interested in $m \ge 2$, since we already know $N_{n,1}(V) = N_n(V)$.

Let $G = \mZ/{q_1}\mZ \times \cdots \times \mZ/{q_m}\mZ$ be a finite Abelian group with $m$ invariant factors: $q_{m} > 1$ and each $q_{i+1}$ divides $q_i$.
By analogy with Section~\ref{sec:proof1}, we say that two $n$-dimensional 
vectors $\vec{a}, \vec{b} \in G^n$ are {\it equivalent\/} modulo $G$,
if  there is an automorphism $\tau$ of $G$ such that 
$$
b_i = \tau(a_i), \qquad 1 \le i \le n.
$$
We also say that a vector $\vec{a}=(a_1,\ldots, a_n) \in G^n$ 
is {\it primitive\/} modulo $G$ if the components $a_1,\ldots, a_n$ generate $G$. 

Let $A_{n}(G)$ be the number of distinct non-equivalent primitive modulo $G$
vectors $\vec{a}=(a_1,\ldots, a_n) \in G^n$. 
Then:
 $$N_{n,m}(V) = \sum_{G\in \fG_m(V)} A_{n}(G).$$
 Indeed,  $L \in \cI_{n,\#G}$ satisfies $\mZ^n/L \simeq G$ if and only 
 there exists a vector  $\vec{a}=(a_1,\ldots, a_n) \in G^n$   primitive modulo $G$,
such that:
$$
L  = \{\vec{x}=(x_1,\ldots, x_n) \in \Z^n~:~
a_1 x_1 + \ldots + a_nx_n = 0 \ \mathrm{in}\ G \}. 
$$
 Note that when $n$ is sufficiently large with respect to $\#G$,  most elements of $G^n$ are primitive modulo $G$: more precisely, 
 Pak~\cite{Pa99} shows that for any $k > 0$, 
 if $(g_1,\dots,g_n) \in G^n$ is picked uniformly at random
 where $n > (k+1) \log \#G +2$,
 then  $g_1,\dots,g_n$ generate the whole group $G$ with probability at least 
 $1-1/\#G^k$,
 which implies that:
 $$A_{n}(G) = \frac{\#G^n ( 1 + O(1/\#G^k))}{\#\mathrm{Aut}(G)}.$$
 In particular, 
 $$A_{n}(G) \sim \frac{\#G^n}{\#\mathrm{Aut}(G)}$$
 as $\# G \to \infty$. 
 Now, there are classical formulas for $\#\mathrm{Aut}(G)$ when $G$ is a finite Abelian group  (see~\cite{HiRh07,Le09,Le10}):
 \begin{fact}
 \label{frac:G p-group}
 If 
 $$
 G = \prod_{i=1}^k (\mZ/{p^{e_i}}\mZ)^{r_i}
$$ 
is a finite Abelian $p$-group in standard form,
that is, $k \ge 0$, $e_1 > \cdots > e_k > 0$, $r_i > 0$,
then:
$$  \#\mathrm{Aut}(G) =  \left( \prod_{i=1}^k \left( \prod_{s=1}^{r_i} (1-p^{-s}) \right)  \right) \left( \prod_{1 \le i,j \le k} p^{\min(e_i,e_j) r_i r_j} \right).$$
 \end{fact}

 \begin{fact}
 \label{frac:G1 G2} If $G_1$  (resp. $G_2$) is a finite Abelian $p_1$-group (resp. $p_2$-group),
  with distinct primes $p_1 \ne p_2$,
then $$  \#\mathrm{Aut}(G_1 \times G_2) =   \#\mathrm{Aut}(G_1)  \times \#\mathrm{Aut}(G_2).$$
\end{fact}

Thus, given a decomposition of $G$ as a product of $p$-groups,
we have a closed formula for $ \#\mathrm{Aut}(G)$,
but if one is given an invariant factor decomposition instead,
then one must first convert it into a $p$-groups decomposition.
We conclude with posing several open problems which might be 
addressed using these formulas:

\begin{prob}
\label{prob:NnmV}
Obtain an  asymptotical formula for $N_{n,m}(V)$ for $m \ge 2$.
\end{prob}

\begin{prob}  Show that  $N_{n,m}(V) = o( \#I_{n,\le V}) $ when $m$ is sufficiently large.
\end{prob}

 \begin{prob}
 \label{prob:Ajt} Show that 
the Ajtai~\cite{Aj96} lattices, that are the classes $\cL_{n,G}$ for $G = (\mZ/q\mZ)^m$,
 form a negligible fraction of $I_{n,\le V}$.
\end{prob}

Here we make a few comments regarding Problem~\ref{prob:Ajt}.
For $G = (\mZ/q\mZ)^m$, where 
$$
q =\prod_{i=1}^\nu  p_i^{e_i}
$$
is the prime number factorisation of $q$, 
each $p$-group of $G$ is of the form $(\mZ/{p_i^{e_i}}\mZ)^{m}$
whose automorphism group, by Fact~\ref{frac:G p-group} has order 
$$
\#\mathrm{Aut}((\mZ/p_i^{e_i}\mZ)^{m}) 
=  p_i^{e_i m^2} \prod_{s=1}^{m} (1-p_i^{-s}), \qquad i =1, \ldots, \nu.
$$
Hence,  
then by Fact~\ref{frac:G1 G2} 
$$\#\mathrm{Aut}((\mZ/q\mZ)^m) =q^{m^2} \prod_{i=1}^\nu  \prod_{s=1}^{m} (1-p_i^{-s}).$$
Note that $\mathrm{Aut}((\mZ/q\mZ)^m)  \simeq \mathrm{GL}_m (\mZ/q\mZ)$.
We also note that due to the application in worst-case to average-case reductions
(see~\cite{Aj96}) we are especially interested in the regime where $n$ is of 
order $m \log m$ and $q$ is of order $m$.

Note that~\cite{GIN} shows how to efficiently sample a random lattice 
from the uniform distribution over $\cL_{n,G}$, given any $G$ for which the factorization of $\#G$ is known.
The knowledge of approximate values of $N_{n,m}(V)$,
see Problem~\ref{prob:NnmV}, 
may lead to an effective sampling of 
 random lattice 
from the uniform distribution over $\cI_{n,\le V}$,
which is an open problem.

\ignore{
It is natural to ask about the number $N_{n,m}(V)$ of 
distinct lattices defines by a system of $m$ congruences, that is,
of lattices
\begin{equation*}
\begin{split}
\cL_n(q, \vec{A} ) = \{\vec{x} & =  ( x_1,   \ldots, x_n) \in \Z^n~:\\
&~a_{1,i} x_1 + \ldots + a_{n,i}x_n  \equiv 0  \pmod {q_i}, \ i =1, \ldots, m\},
\end{split}
\end{equation*}
defined an $m\times n$ integer matrix $\vec{A}= (a_{i,j})_{i,j=1}^{m,n} \in \Z^{m\times n}$.

Furthermore,  it is also interesting to count 
the number of $n$-dimensional lattices of rank $r$ with a
squarefree determinant. We note that the ideas of Poonen~\cite{Poon} may 
eventually lead to counting squarefree determinants
of integral matrices, ordered by the size of their 
entries. However for  algorithmic applications it is more relevant 
to count squarefree determinants
of integral matrices ordered by the size of the  determinants,
which seems to be a very difficult question. 
}

\ignore{
\subsection{Comments (food for thought)}

I don't know $\# \fG_m(V)$ except obviously $\# \fG_1(V)$, but here is what is known.
Let $a(n)$ denote the number of non-isomorphic Abelian groups of order $n$.
It is known that
$$  \sum_{m \ge 1} \# \fG_m(x) = \sum_{n \le x} a(n) = A_1 x + A_2 x^{1/2} + A_3 x^{1/3} + R(x),$$
where $A_j = \prod_{k=1}^{\infty} \zeta(k/j)$ ($j=1,2,3$),
and the best result for the error term is $R(x) \ll x^{1/4+\varepsilon}$ for any $\varepsilon > 0$
(see~\cite{RoSa06}).
We have  $A_1 \approx 2.29\dots$.
}

\section*{Acknowledgements}

P.~Q.~Nguyen would like to thank Christophe Delaunay for very helpful discussions on~\cite{CoLe84}.

During the preparation of this paper, 
P.~Q.~Nguyen was supported in part by  China's 973 Program, Grant 2013CB834205,  and NSFC's Key Project, Grant 61133013,
and I.~E.~Shparlinski was supported in part by   ARC grants DP130100237 and~DP140100118.


\begin{thebibliography}{99}

\bibitem{Aj96}
M.~Ajtai, `Generating hard instances of lattice problems', 
{\it  Proc.   28th  ACM Symp. on Theory of Comp.\/}, ACM, 1996, 99--108. 

\bibitem{CoLe84}
H.~Cohen and H.~W. Lenstra, Jr., 
`Heuristics on class groups of number fields',
{\it Lect. Notes in Math\/}, Springer-Verlag,
Berlin, {\bf 1068} (1984), 33--62.  

\bibitem{FuGo15}
J.~Fulman and L.~Goldstein,
`Stein's method and the rank distribution of random matrices over
  finite fields', {\it  Ann. Probab.\/}, {\bf 43} (2015),  1274--1314. 



\bibitem{GIN} N. Gama,  M. Izabachene, P.~Q. Nguyen and X. Xie,
`Structural lattice reduction: Generalized worst-case to average-case 
reductions and homomorphic cryptosystems',  
{\it Cryptology ePrint Archive: Report 2014/283\/},  
2014
(available from \url{http://eprint.iacr.org/2014/283}). 

\bibitem{Ge08}
E.-U. Gekeler, `Statistics about elliptic curves over finite prime fields', 
{\it Manuscr. Math.\/}, {\bf 127} (2008), 55--67.



\bibitem{GiGr06}
H.~Gillet and D.~R. Grayson, 
`Volumes of symmetric spaces via lattice points', 
{\it Doc. Math.\/}, {\bf 11}  (2006),  425--447.

\bibitem{HW} G. H. Hardy and E. M. Wright, {\it An introduction
to the theory of numbers\/}, Oxford Univ. Press, Oxford, 1979.

\bibitem{HiRh07}
C.~J. Hillar and D.~L. Rhea, 
`Automorphisms of finite abelian groups', 
{\it Amer. Math. Monthly\/}, {\bf 114} (2007), 917--923.

\bibitem{Hool} C. Hooley, 
`A note on square-free numbers in arithmetic progressions', 
{\it Bull. London Math. Soc.\/}, {\bf 7} (1975), 133--138.

\bibitem{Lan} E. Landau, `\"{U}ber die Zahlentheoretische
Function $\varphi(n)$ und ihre Beziehung zum Goldbachschen Satz',
{\it Nachr. K\"{o}niglichen Ges. Wiss. G\"{o}ttingen, Math.-Phys.
Klasse\/}, G\"{o}ttingen, 1900, 177--186.

\bibitem{Le09}
J.~Lengler, {\it The {C}ohen-{L}enstra Heuristic for Finite Abelian Groups\/},
PhD thesis, Universit\"at des Saarlandes, 2009.


\bibitem{Le10}
J.~Lengler, `The Cohen-Lenstra heuristic: Methodology and results',
{\it J. Algebra\/}, {\bf 323} (2010), 2960--2976.

\bibitem{Mont} H. L. Montgomery,
`Primes in  arithmetic progressions', 
{\it Mich. Math. J.\/}, {\bf 17} (1970), 33--39.

\bibitem{Pa99}
I.~Pak, `On probability of generating a finite group', 
{\it Preprint\/}, 1999.


\bibitem{PaSch} A. Paz and C.-P. Schnorr,
`Approximating integer lattices by lattices with cyclic factor groups', 
{\it Lect. Notes in Comp. Sci.\/}, Springer-Verlag,
Berlin, {\bf 267} (1987), 386--393.


\bibitem{Prach} K. Prachar, `\"Uber die kleinste quadratfreie Zahl einer arithmetischen Reihe', {\it Monatsh. Math.\/}, {\bf 62} (1958), 173--176.

\bibitem{RoSa06}
O.~Robert and P.~Sargos, 
`Three-dimensional exponential sums with monomials', 
{\it J. Reine Angew. Math.\/}, {\bf 591} (2006), 1--20.



\bibitem{Sc68}
W.~M. Schmidt, 
`Asymptotic formulae for point lattices of bounded determinant and
  subspaces of bounded height', 
{\it Duke Math. J.\/}, {\bf 35} (1968), 327--339.



\bibitem{Trol} M. Trolin, `Lattices with many cycles are dense', 
{\it Lect. Notes in Comp. Sci.\/}, Springer-Verlag,
Berlin, {\bf 2996} (2004),  370--381.

\bibitem{Wa27}
D.~R. Ward, `Some series involving Euler's function',
{\it J. London Math. Soc.\/}, {\bf 210} (1927), 1--2.



\end{thebibliography}
\end{document}